\documentclass[final, 3p, times, authoryear]{elsarticle}
\usepackage{amssymb}

\usepackage{amsthm}

 \usepackage{lineno}

\usepackage{setspace}
\usepackage{float}

\biboptions{sort}

\usepackage{amsmath, amsthm, amssymb} % American Mathematical Society packages
\usepackage{bm} % For bold mathematic notation

\usepackage{epsfig} % For figures
\usepackage{graphicx} % For graphics

\def\bea{\begin{eqnarray}}
\def\eea{\end{eqnarray}}
\def\bean{\begin{eqnarray*}}
\def\eean{\end{eqnarray*}}
\renewcommand\eqref[1]{(\ref{#1})}

\def\N{\mathbb{N}}
\def\Z{\mathbb{Z}}

\def\Ind#1{\,\mathbb{I}\{#1\}\,}

\theoremstyle{plain}
\newtheorem{theorem}{Theorem}
\newtheorem{proposition}[theorem]{Proposition}

\theoremstyle{definition}

\theoremstyle{remark}

\theoremstyle{note}

\usepackage{comment}
\usepackage{ulem} % used for \sout
\usepackage[usenames,dvipsnames]{xcolor}%\usepackage{color}
\definecolor{orange}{rgb}{1,0.5,0}
\definecolor{green}{rgb}{0.513,0.73,0.442}
\renewcommand\emph[1]{{\it #1}} % the ulem package redefines \emph
\ifx \q \undefined
  \def \q[#1][#2][#3]{q_{#2}^{#3}\left(#1\right)}% for std generator
\def \figref#1{\ref{#1}}% for figure reference

\def\N{\mathbb{N}} %% natural numbers
\def\Z{\mathbb{Z}} %% integer numbers
\def\dt{{h}} % for time increment
\def\Ind#1{\,\mathbb{I}\{#1\}\,}  %% indicator function
\def\dt{{h}} % for time increment
\def\dtt{{\bar{h}}} % for bounding interval
\def\MultInd{{c}}
\ifx \q \undefined
  \def \q[#1][#2][#3]{q_{#2}^{#3}\left(#1\right)}% for std generator
  \def \qt[#1][#2][#3]{\tilde q_{#2}^{#3}\left(#1\right)}% for tilde generator
\fi

\newcommand{\GammaFunc}{{G}}

\newcommand{\dquery}[1]{{\relax}}
\newcommand{\daihai}[1]{{\relax}}

\def \figref[#1]{\ref{#1}}% for figure reference

\def\recovery{r}
\def\Io{I_2}
\def\Ii{I_1}
\def\Iovar{{I^*_{2}}}
\def\Iivar{{I^*_{1}}}
\def\noise{\xi}
\def\Sb{S}
\def\So{S_2}
\def\Si{S_1}
\def\i{1}
\def\o{2}
\def\reservoir{\omega}
\def\csep{ }  %separation, e.g. \mu_{S \csep Io}

\excludecomment{repeatIntro}% NO
\excludecomment{styleRevisionFormat}% NO
\excludecomment{styleRevisionColor}% NO
\definecolor{NeonYellowWhiteOnBlack}{rgb}{0.016,0.009, 0.975}
\def\kwy#1{#1}
\begin{styleRevisionColor}
\def\kwy#1{\textcolor{NeonYellowWhiteOnBlack}{#1}}
\end{styleRevisionColor}

\definecolor{NeonBlueWhiteOnBlack}{rgb}{0.975, 0.016, 0.009}
\def\kwb#1{#1}
\begin{styleRevisionColor}
\def\kwb#1{\textcolor{NeonBlueWhiteOnBlack}{#1}}
\end{styleRevisionColor}

\definecolor{NeonOrangeWhiteOnBlack}{rgb}{0.0075, 0.44, 0.5525}
\def\kwo#1{#1}
\begin{styleRevisionColor}
\def\kwo#1{\textcolor{NeonOrangeWhiteOnBlack}{#1}}
\end{styleRevisionColor}

\def\kwdb#1{#1}
\begin{styleRevisionColor}
\def\kwdb#1{\textcolor{orange}{#1}}
\end{styleRevisionColor}

\def\kwp#1{#1}
\begin{styleRevisionColor}
\def\kwp#1{\textcolor{OliveGreen}{#1}}
\end{styleRevisionColor}

\definecolor{reproduction-gray}{gray}{0.65}

\begin{styleRevisionFormat}

\end{styleRevisionFormat}
\definecolor{summarizing-gray}{gray}{0.35}

\usepackage[geometry]{ifsym}
\newcommand{\Ring}[1]{\raisebox{-1pt}{\begin{tabular}{@{}c@{}}{\small #1}\\
      [-11.5pt]\BigCircle\end{tabular}}}
\newcounter{sectionSentenceCounter}[section] % to be increased with each new summarizing sentence in section introductory paragraphs (e.g., \)
\newcounter{subsectionSentenceCounter}[subsection] 
\newcommand{\secsentence}{\refstepcounter{sectionSentenceCounter}\Ring\thesectionSentenceCounter}
\begin{document}

\begin{frontmatter}
%\title{Co-jumps and infinitesimal covariances of Markov counting systems}
%\title{Co-jumps in population dynamics subject to random environments}
%\title{Co-jumps and population dynamics in random environments}
\title{Co-jumps and Markov counting systems in random environments}
%front matter mark
%
\author{Carles Bret\'{o}} \fnref{label1}
\fntext[label1]{Tel. +34916245855; Fax:+34916249848}
\ead{carles.breto@uc3m.es}
%\address{Departamento de Estad\'{i}stica and Instituto Flores de Lemus, Universidad Carlos III de Madrid, C/ Madrid 126, Getafe, 28903, Madrid, Spain}
\address{Departamento de Estadistica and Instituto Flores de Lemus, Universidad Carlos III de Madrid, C/ Madrid 126, Getafe, 28903, Madrid, Spain}
%
%\author{\paperauthor}
%\address{}
%
\begin{abstract}  
We provide transition rates for Markov counting systems subject to correlated environmental noises motivated by multi-strain disease models. 
Such noises induce simultaneous counts, which can help model infinitesimal count correlation (regardless of whether such correlation is due to correlated noises). %(regardless of whether it comes from correlated noises or not).
%40 words
%
%We provide transition rates for Markov counting systems subject to correlated environmental noises motivated by multi-strain disease models. 
%Such noises induce simultaneous infections of individuals, which we show can help model count infinitesimal correlation (due or not to correlated noises).
%40 words
%
%We provide closed-form transition rates for Markov counting systems subject to correlated environmental noises motivated by multi-strain epidemiological models. 
%Such noises induce simultaneous counts (co-jumps), which we show in turn induce (and are necessary for) infinitesimal correlations between counting processes.
%40 words
%
%We show that introducing correlated noises to the rates of several Markov counting processes induces simultaneous counts (or co-jumps), which in turn induce (and are necessary for) infinitesimal correlation between such processes, as we illustrate with multi-strain epidemiological models.
%39 words
%
%We show that introducing correlated noises to the rates of interacting Markov counting processes induces simultaneous counts (or co-jumps), which in turn induce (and are needed for) infinitesimal correlation between such processes and which we illustrate with multi-strain epidemiological models.
%40 words
%
%We show that co-jumps are necessary and sufficient for interacting Markov counting processes to be infinitesimally correlated. 
%This occurs for example when the rates of such processes are affected by correlated environmental noises, as we illustrate with multi-strain epidemiological models.
%40 words
%
\end{abstract}
\end{frontmatter}
%%\linenumbers

%BODY
%\section{}
%\label{}
%
%%%%%%%%%%%%%%%%%%%%%%%%%%%%%%%%%%%%%%%%%%%%%%%%%%%%
%%%%%%%%%%%%%%%%%%%%%%%%%%%%%%%%%%%%%%%%%%%%%%%%%%%%
%
% introduction mark
\section{Introduction}
\label{sec:intro}
%SEEMING TRUTH
\def \SeemingTruth{\kwp{Continuous-time stochastic processes} have proved to be useful for research in many areas of science. 
Such 
\kwp{processes} are naturally defined by 
\kwdb{infinitesimal parameter functions}, like transition rates in the case of Markov chains. %or infinitesimal mean and covariance functions (for diffusion processes).
Assuming that such 
\kwdb{parameter functions} are subject to external noise has been referred to as a random environment and has proved to be a useful approach to \kwp{data analysis}. 
}%SEEMING TRUTH DETAILS
\def \SeemingTruthDetails{A recent example of such 
\kwp{applied work} is that in \cite{shrestha2011}. 
This work studied the dynamics governing interactions among multiple infectious pathogens; 
it found promising results based on Markov counting systems for which multiple transition rates are subject to a single, common \kwdb{external noise}. 
Such common 
\kwdb{noise} results in correlation between these transition rates. 
Subjecting transition rates to such noise is known to result in a new Markov counting system defined by
\kwp{new rates}. 
}%PROBLEM
\def\ConditionOne{However, these 
\kwp{new rates} have been derived in closed form only for the case of 
\kwdb{independent noises}.
}%SO WHAT?
\def\ConsequenceOne{This lack of closed-form rates for the case of 
\kwdb{correlated noises} creates uncertainty about the correlated system properties and might make applied researchers 
\kwp{reluctant} to take advantage of the promising approach presented by \cite{shrestha2011}.
}%SOLUTION
\def\SolutionOne{To prevent such 
\kwp{reticence}, this paper considers introducing correlated external noises to the rates of Markov counting systems and provides closed-form expressions for 
\kwdb{new rates} that capture the effect of correlated noise.
These 
\kwdb{new rates} are the main contribution of the paper and are based on novel closed-form expressions for the 
\kwp{infinitesimal covariances} of Markov counting systems. 
In addition, this unusual focus on 
\kwp{system covariances} provides an alternative interpretation of correlated noises in terms of simultaneous transitions in non-random environments. 
}
% PUT ALL THE TEXT ABOVE IN THE CORRECT STRUCTURE
%
%
%SEEMING TRUTH
%
%
\begin{styleRevisionFormat}
\def \seemingTruth{[Seeming truth:]}
\kwb{\seemingTruth}
\end{styleRevisionFormat}
\begin{styleRevisionFormat}% begin comment environment to include/exclude style revision format
\secsentence\label{SeemingTruth}
\end{styleRevisionFormat}
\SeemingTruth
%
%
%SEEMING TRUTH DETAILS
%
%
\begin{styleRevisionFormat}
\def \seemingTruthDetails{[Seeming truth details:]}
\kwb{\seemingTruthDetails}
\end{styleRevisionFormat}
\begin{styleRevisionFormat}% begin comment environment to include/exclude style revision format
\secsentence\label{SeemingTruthDetails}
\end{styleRevisionFormat}
\SeemingTruthDetails
%
%
%PROBLEM: STATE A CONDITION
%
\begin{styleRevisionFormat}
\def \problem{[Problem condition:]}
\kwb{\problem}
\end{styleRevisionFormat}
\begin{styleRevisionFormat}% begin comment environment to include/exclude style revision format
\secsentence\label{ConditionOne} 
\end{styleRevisionFormat}
\ConditionOne
%
%
% PROBLEM CONSEQUENCE: IMAGINARY ``SO WHAT?''
%
%
\begin{styleRevisionFormat}
\def \soWhat{[So what?]}
\kwb{\soWhat}
\end{styleRevisionFormat}
\begin{styleRevisionFormat}% begin comment environment to include/exclude style revision format
\secsentence\label{ConsequenceOne}
\end{styleRevisionFormat}
\ConsequenceOne
%
%
% PROBLEM SOLUTION
%
%
\begin{styleRevisionFormat}
\def \solution{[Solution:]}
\kwb{\solution}
\end{styleRevisionFormat}
\begin{styleRevisionFormat}% begin comment environment to include/exclude style revision format
\secsentence\label{SolutionOne}
\end{styleRevisionFormat}
\SolutionOne
%
%
%%%%%%%%%%%%%%%%%%%%%%%%%%%%%%%%%%%%%%%%%%%%%%%%%%%%%%%%%%%%%%%%
% BEGIN INTRODUCTION PARAGRAPHS THAT DEVELOP THE IDEAS SUMMARIZED BY THE FIRST PARAGRAPH
%%%%%%%%%%%%%%%%%%%%%%%%%%%%%%%%%%%%%%%%%%%%%%%%%%%%%%%%%%%%%%%%
%
¥% extended introduction mark
%\\% INTEPARAGRAPH SPACE 
%
%\begin{comment}
\begin{repeatIntro}
\\% INTEPARAGRAPH SPACE 
\end{repeatIntro}

%INTEPARAGRAPH SPACE
\begin{repeatIntro}% begin comment environment to include/exclude repeated introduction sentences before its corresponding paragraph
\color{reproduction-gray} % set the text color to "reproduction-gray" to point out that the text is reproducing a sentence from the first paragraph of the introduction (which is supposed to summarize the whole introduction section).
\begin{styleRevisionFormat}% begin comment environment to include/exclude style revision format
\Ring{\ref{SeemingTruth}} \seemingTruth
\end{styleRevisionFormat}
\SeemingTruth \\
\\ % break line so to separate reproduced text from new text (that remains in the same color)
\color{black} % return to black color at the end of the block in "reproduction color" that is formed by the reproduced sentence from the first paragraph and  by the current paragraph point sentence.
\end{repeatIntro}
¥% seeming truth paragraph mark 
%\paragraphTopicSentence{\indent 
The study of properties of \kwp{counting processes} has benefited among others the fields of \kwy{epidemiology and ecology}, which have relied on \kwo{Markov counting systems} both in deterministic and \kwb{stochastic} environments (the latter sometimes being favoured by \kwp{empirical evidence}). 
%} 
%
% REST OF THE PARAGRAPH
Research in these two \kwy{disciplines} has taken advantage of theoretical investigation of counting systems both historically \citep{kermack1927, bartlett1956} %, pharmacokinetics \citep{%matis1979, srivastava2002} and engineering and operations research \citep{%doig1957, jackson2002}. 
and more recently \citep{cauchemez2008, he2010}.
Such counting systems can often be seen as systems of interacting Markov counting processes or Markov counting systems \citep{breto2011}, which include %which range from plain birth-death processes to more complex processes, %(that may consist of multiple interacting birth-death processes, each with its own sub-population), like 
networks of queues \citep{bremaud1999} and compartmental models \citep{jacquez1996, matis2000}. 
%\kwo{Markov counting systems} are very flexible: not only can they model concurrent birth (or death) events of several individuals within a given sub-population (i.e., allow for compound events) but also allow for simultaneous events across different sub-populations (i.e., allow for co-jumps). 
%Nevertheless, most systems considered in the literature rule out both compound events and co-jumps, which conveniently simplifies the \kwdb{transition rates} of the Makov counting system (as we precise in Section~\ref{sec:rates}). 
\kwo{Markov counting systems} are Markov chains and are hence naturally defined by transition rates. 
Noisy transition rates are often referred to as environmental stochasticity in epidemiology and ecology \citep{engen1998}.  
The role of \kwb{such stochasticity} has been extensively studied, including in the context of deterministic ODE skeletons driven by diffusions \citep{dureau2013, hu2011, ionides2006-pnas, king2008} and driven by L\'{e}vy processes \citep{bhadra2011, laneri2010}. 
The role of stochastic environments has also been studied in the context of Markov counting systems, both paying attention to the system probabilistic properties \citep[e.g.,][]{marion2000, varughese2008, breto2009} and focusing on the biological implications for \kwp{applications} \citep[e.g.,][]{shrestha2011}.
%
% END OF PARAGRAPH
%
\begin{repeatIntro}
\\% INTEPARAGRAPH SPACE 
\end{repeatIntro}

%INTEPARAGRAPH SPACE
\begin{repeatIntro}% begin comment environment to include/exclude repeated introduction sentences before its corresponding paragraph
\color{reproduction-gray} % set the text color to "reproduction-gray" to point out that the text is reproducing a sentence from the first paragraph of the introduction (which is supposed to summarize the whole introduction section).
\begin{styleRevisionFormat}% begin comment environment to include/exclude style revision format
\Ring{\ref{SeemingTruthDetails}} \seemingTruthDetails
\end{styleRevisionFormat}
\SeemingTruthDetails \\
\\ % break line so to separate reproduced text from new text (that remains in the same color)
\color{black} % return to black color at the end of the block in "reproduction color" that is formed by the reproduced sentence from the first paragraph and  by the current paragraph point sentence.
\end{repeatIntro}
¥% seeming truth paragraph mark 
%\paragraphTopicSentence{\indent 
Epidemiological \kwp{applications} have come to consider \kwy{multiple interacting pathogens} and to study them based on counting systems subject to genuinely \kwo{correlated environmental} noise, fitting in the framework provided by \kwb{\citet{breto2009}} who formalize the \kwp{transition rates} of the system subject to noise.
%}
%
% REST OF THE PARAGRAPH
\kwy{Pathogen interaction} has received attention for some time now \citep{fenton2005, kamosasaki2002}, both without considering the role of external noises \citep{aguiar2011, buckee2011, reich2013} and considering it \citep{breto2009, shrestha2011}. 
In particular, \citet{shrestha2011} consider a Markov counting system corresponding to a compartmental model of the susceptible-infectious-recovered type (a simpler version of which is considered by \cite{breto2009} and reproduced in Figure~\ref{fig:model}). 
In \cite{shrestha2011}, two pathogens co-exist but there are more than two possible different types of infection (depending on the history of past infections of individuals). 
The rate at which these different types of infection occur are assumed to be subject to a single (common to all infection types) external white noise, making all infection rates correlated. 
Such rate \kwo{correlation} has been formalized by \kwb{\citet{breto2009}}, showing that the system subject to noises is a new Markov counting system not only when the noises are independent (in which case they even provide closed-form rates) but also when the noises are correlated, on which this paper \kwdb{focuses}. 
%
% END OF PARAGRAPH
%
\begin{repeatIntro}
\\% INTEPARAGRAPH SPACE 
\end{repeatIntro}
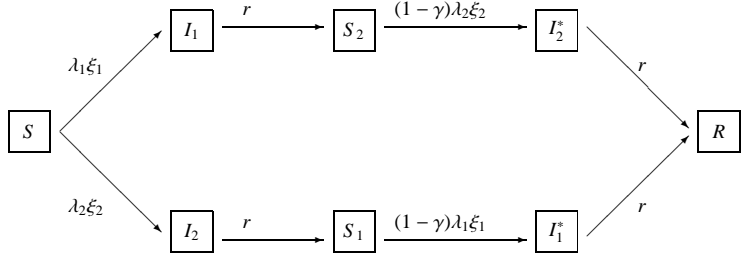
\begin{figure}
\def\Clabel{C} % label for compartment membership function
\def\Ilabel{\zeta} % label for individuals in compartment model

\def\typeA{{\prime}}
\def\typeB{{\prime\prime}}
\def\typeC{{\prime\prime\prime}}
\def\matchSet{S} % set on which two approximations to N_ij(t) agree

\def\recovery{r}
\def\Io{I_2}
\def\Ii{I_1}
\def\Iovar{{I^*_{2}}}
\def\Iivar{{I^*_{1}}}
\def\noise{\xi}
\def\Sb{S}
\def\So{S_2}
\def\Si{S_1}
\def\i{1}
\def\o{2}
\def\reservoir{\omega}
\def\csep{ }  %separation, e.g. \mu_{S \csep Io}

\def\envstoch{\sigma}
\def\demstoch{\phi}
\def\preport{\rho}
\def\Ci{C_{1,t}}
\def\Co{C_{2,t}}

\def\GammaFunc{G}

        \resizebox{10cm}{!}
            {
            \begin{picture}(370,115)(-70,10)
            %BIRTH
%                \put(-70,50){\framebox(20,20){$B$}}
%                \put(-45,60){\vector(1,0){50}} %arrow B-Xn
%                \put(-35,66){$b(t)$} % rate in arrow B-Xn
            %S_b
                \put(10,50){\framebox(20,20){$\Sb$}}
                \put(35,60){\vector(1,1){50}} %arrow up Xn
                \put(39,90){$\lambda_{\i}\noise_{\i}$}
                \put(35,60){\vector(1,-1){50}} %arrow down Xn
                \put(39,20){$\lambda_{\o}\noise_{\o}$}
%                \color{gray}\put(35,60){\vector(1,0){100}} %arrow Xn-D
%                \put(55,66){$\mu$}\color{black} % rate in arrow Xn-D
            %I_b
                \put(90,100){\framebox(20,20){$\Ii$}}
                \put(115,112){\vector(1,0){50}} %arrow Ib-Sb
                \put(125,118){$\recovery $}
            %I_g
                \put(90,0){\framebox(20,20){$\Io$}}
                \put(115,6){\vector(1,0){50}} %arrow Ig-Sg
                \put(125,12){$\recovery $}
%                \put(142,90){$(1-\gamma)\lambda_{\o}$}
%                \put(142,25){$(1-\gamma)\lambda_{\i} $}
%                \put(100,95){\vector(0,-1){20}} %arrow Xi-D
%                \put(107,87){$\mu $} % rate of arrow Xi-D
%                \put(100,25){\vector(0,1){20}} % arrow Xo-D
%                \put(107,31){$\mu $} % rate of arrow Xo-D
%                \put(90,50){\framebox(20,20){$D$}}
            %S_o
                \put(170,100){\framebox(20,20){$\So$}}
                \put(195,112){\vector(1,0){70}} %arrow Sb-Ig'
                \put(200,118){$(1-\gamma)\lambda_{\o}\noise_{\o}$}
%                  \put(195,98){\vector(4,-1){145}}
%                  \put(215,72){$\nu\gamma\lambda_{\o}\noise_{\o}$}
            %S_i
                \put(170,0){\framebox(20,20){$\Si$}}
                \put(195,6){\vector(1,0){70}} %arrow Sg-Ig'
                \put(200,12){$(1-\gamma)\lambda_{\i}\noise_{\i}$}
 %                 \put(195,22){\vector(4,1){145}}
 %                 \put(215,42){$\nu\gamma\lambda_{\i}\noise_{\i}$}
%                \put(165,60){\vector(-1,0){50}} %arrow  Xb-D
%                \put(125,66){$\mu $} % rate in arrow  Xb-D
            %I_o'
                \put(270,100){\framebox(20,20){$\Iovar$}}
                \put(295,112){\vector(1,-1){50}} %arrow Ig'-Rb
                \put(320,90){$\recovery $}
            %I_i'
                \put(270,0){\framebox(20,20){$\Iivar$}}
                \put(295,8){\vector(1,1){50}} %arrow Ig'-Rb
                \put(320,20){$\recovery$}
            %R_b
                \put(350,50){\framebox(20,20){$R$}}
            \end{picture}
        }%end resizebox
\caption{Multi-strain SIR-type compartmental model of \cite{breto2009}. 
This model will be used in Section~\ref{sec:application} to illustrate our results. 
%Flow diagram for cholera, including interactions between the two major serotypes. 
Each individual falls in one compartment:
$\Sb$, susceptible to both strains %Inaba and Ogawa serotypes
; $\Ii$, infected with strain 1%Inaba
; $\Io$, infected with strain 2%Ogawa
; $\Si$, susceptible to strain 1 %Inaba 
(but immune to strain 2%Ogawa
); $\So$, susceptible to strain 2 %Ogawa 
(but immune to strain 1%Inaba
); $\Iivar$, infected with strain 1 %Inaba 
(but immune to strain 2%Ogawa
); $\Iovar$, infected with strain 2 %Ogawa 
(but immune to strain 1%Inaba
); and $R$, immune to both strains%serotypes
. 
Regarding demography, births enter $S$ from compartment $B$ (not plotted), at rate $b(t)$ driven by birth data (which is treated as a covariate), and all individuals have a common mortality rate $m$ at which they leave each compartment in the diagram into $D$ (not plotted). 
Regarding disease dynamics, $r$ is the recovery rate from infection; $\gamma$ measures the strength of cross-immunity between strains; and $\lambda_{i}$ is the per-capita infection rate of strain $i$ with $\noise_{i}$ being the stochastic noise on this rate. 
Moreover, $\lambda_{i} = \beta(t)( I_{i}(t) + I_{i}^{*}(t) )^\alpha/P(t)+\reservoir$, where $0 \le \beta(t)$ is parameterized with a trend and a smooth seasonal component, $0 \le \reservoir$ models infections from an environmental reservoir and $0 \le \alpha \le 1$ captures inhomogeneous mixing of the population.
}\label{fig:model}
\end{figure}
%
%INTEPARAGRAPH SPACE
\begin{repeatIntro}% begin comment environment to include/exclude repeated introduction sentences before its corresponding paragraph
\color{reproduction-gray} % set the text color to "reproduction-gray" to point out that the text is reproducing a sentence from the first paragraph of the introduction (which is supposed to summarize the whole introduction section).
\begin{styleRevisionFormat}% begin comment environment to include/exclude style revision format
\Ring{\ref{ConditionOne}} \problem
\end{styleRevisionFormat}
\ConditionOne \\
\\ % break line so to separate reproduced text from new text (that remains in the same color)
\color{black} % return to black color at the end of the block in "reproduction color" that is formed by the reproduced sentence from the first paragraph and  by the current paragraph point sentence.
\end{repeatIntro}
¥% seeming truth paragraph mark 
%\paragraphTopicSentence{\indent 
The \kwdb{problem} we take up in this paper is providing closed-form transition rates that define Markov counting systems accounting for correlated noises and is made difficult by the lack of closed-form \kwy{transition probabilities} of \kwo{general systems}, working \kwp{against the inclusion of biologically genuine noise correlation}.
%}
%
% REST OF THE PARAGRAPH
%Transition rates of Markov counting systems can be derived as certain limits of transition probabilities \citep{bremaud1999}. 
Closed-form \kwy{transition probabilities} are readily available for basic systems, e.g., those of a Poisson process correspond to a Poisson distribution %; those of a linear pure birth process, to a negative binomial distribution; 
and those of a linear pure death process to a binomial distribution \citep{bharucha1960}. 
%Now, consider combining this Poisson and binomial processes in to a birth-death process where births occur according to the Poisson process and deaths according to the binomial process
However, they are not available for general compartmental models, including the system of interacting birth-death processes considered by \cite{shrestha2011} or that represented in Figure~\ref{fig:model}. 
If \kwo{such closed-form general system transition probabilities} were available, then the desired closed-form transition rates might be pursued by direct integration of the noise from those (now randomized) probabilities. 
Such direct approach is feasible in basic cases \cite[like the bivariate death process of][]{breto2011} but not in more sophisticated models, where the lack of closed-form rates \kwp{casts a shadow} over the appeal of correlated noises in applications of a realistic degree of complexity. 
%
% END OF PARAGRAPH
%
\begin{repeatIntro}
\\% INTEPARAGRAPH SPACE 
\end{repeatIntro}

%INTEPARAGRAPH SPACE
\begin{repeatIntro}% begin comment environment to include/exclude repeated introduction sentences before its corresponding paragraph
\color{reproduction-gray} % set the text color to "reproduction-gray" to point out that the text is reproducing a sentence from the first paragraph of the introduction (which is supposed to summarize the whole introduction section).
\begin{styleRevisionFormat}% begin comment environment to include/exclude style revision format
\Ring{\ref{ConsequenceOne}} \soWhat
\end{styleRevisionFormat}
\ConsequenceOne \\
\\ % break line so to separate reproduced text from new text (that remains in the same color)
\color{black} % return to black color at the end of the block in "reproduction color" that is formed by the reproduced sentence from the first paragraph and  by the current paragraph point sentence.
\end{repeatIntro}
¥% seeming truth paragraph mark 
%\paragraphTopicSentence{\indent 
A key \kwp{downside} of lacking closed-form rates is that the promising results of \kwy{\cite{shrestha2011}} and the \kwo{genuine biological rationale} behind correlated noises may be outweighted by uncertainty about the \kwb{properties of the model} subject to noise and about the interpretation of empirical results, which we seek to prevent \kwdb{with this paper}.
%}
%
% REST OF THE PARAGRAPH
\kwy{\cite{shrestha2011}} show that it is feasible to arrive at correct and precise biological conclusions regarding pathogen interaction based on their Markov counting systems with correlated noises. 
In addition, a heuristic \kwo{biological justification} for correlated noises could be as follows: while localized environmental variations need not affect all types of infection, changes at a larger scale in the environment should be expected to, like heat or cold waves. 
However, unless the \kwb{properties of the model} after subjecting it to correlated noise are clear and appealing, such correlations might be considered a nuisance or something foreign and hence avoided in actual applications, where empirical findings need to be interpreted (which might be done more confidently in the context of simpler models). 
Providing a tool to help in such interpretation is the ultimate \kwdb{goal of this paper}. 
%
% END OF PARAGRAPH
%
\begin{repeatIntro}
\\% INTEPARAGRAPH SPACE 
\end{repeatIntro}

%INTEPARAGRAPH SPACE
\begin{repeatIntro}% begin comment environment to include/exclude repeated introduction sentences before its corresponding paragraph
\color{reproduction-gray} % set the text color to "reproduction-gray" to point out that the text is reproducing a sentence from the first paragraph of the introduction (which is supposed to summarize the whole introduction section).
\begin{styleRevisionFormat}% begin comment environment to include/exclude style revision format
\Ring{\ref{SolutionOne}} \solution
\end{styleRevisionFormat}
\SolutionOne \\
\\ % break line so to separate reproduced text from new text (that remains in the same color)
\color{black} % return to black color at the end of the block in "reproduction color" that is formed by the reproduced sentence from the first paragraph and  by the current paragraph point sentence.
\end{repeatIntro}
¥% seeming truth paragraph mark 
%\paragraphTopicSentence{\indent 
The \kwdb{main contribution} of this paper is to provide \kwy{closed-form transition rates} for Markov counting systems subject to correlated noises based on the system \kwo{infinitesimal covariances} and to provide an \kwb{illustration} in the context of biological analysis of multi-strain pathogen dynamics. 
%}
%
% REST OF THE PARAGRAPH
%
The provided \kwy{closed-form expressions} %are provided in Section~\ref{sec:application} and 
apply to a broad range of cases considered in the applied literature. 
They reduce the uncertainty about the model properties by giving a precise definition of the system as formalized in Section~\ref{sec:notation}. 
%They also provide a quantification of the role of external noise via infinitesimal moments. % and facilitate simulation-based statistical inference by enabling exact simulation of sample paths. 
In addition, they are motivated by the system infinitesimal covariances derived in Sections~\ref{sec:bivar} and~\ref{sec:inf-cov}, which allow circumventing the above mentioned problem of unavailable transition probabilities from which to directly integrate out the noise. 
Our focus on \kwo{infinitesimal covariances} is unusual in the context of Markov counting systems (although as natural as in the context of multivariate diffusions) and leads to the novel closed-form expressions for them provided in Theorem~\ref{thm:inf-cov}. 
These expressions %provide an alternative interpretation of empirical results favouring random environments. They 
show that correlated noises induce simultaneous counts and that these in turn induce stronger correlations within the system. 
Hence, if additional correlation is demanded by data, it could be modelled with random environments. 
In this case, these environments could be interpreted as devices that generate the needed correlation in a non-random environment, instead of as actual random changes in parameters (very much like parameter randomization can be interpreted as a device to generate over-dispersion).
This is illustrated in Section~\ref{sec:application}, \kwb{where} the rates and interpretation of the role of correlated noises for Figure~\figref[fig:model] are given. 
%\kwb{time randomization} is considered for several processes of interest in the biological and social sciences and applied in detail to the widespread Poisson, linear birth and linear death processes to construct novel, infinitesimally over-dispersed processes and multivariate systems without the interpretation issues of time-changed models defined only implicitly. 

%These 
%\kwdb{new rates} are the main contribution of the paper and shed light on the system properties besides being useful to applied researchers in several ways, including providing an interpretation of correlated noises in terms of infinitesimal covariances between counting processes of the system. 
%
% END OF PARAGRAPH
%
%
%%\\ %INTEPARAGRAPH SPACE
%%
%%%INTEPARAGRAPH SPACE
%
%
%%%%%%%%%%%%%%%%%%%%%%%%%%%%%%%%%%%%%%%%%%%%%%%%%%%%
%%%%%%%%%%%%%%%%%%%%%%%%%%%%%%%%%%%%%%%%%%%%%%%%%%%%
%
\section{Markov counting systems without external noise} 
\label{sec:notation} 
%Let us first \kwp{define} Continuous-time Markov counting systems and give an \kwdb{example} (in Figure~\figref[fig:model]) before considering the \kwp{effect of correlated noises}. % on the \kwdb{infinitesimal covariances} of the system. 
%INTEPARAGRAPH SPACE

%INTEPARAGRAPH SPACE
Markov counting systems are \kwp{defined} as Markov chains driven by a collection of interacting \kwy{counting processes} that fully characterizes the \kwo{transition rates} of the system \citep{breto2009} and such definition can often be formalized in a \kwdb{diagram} (similar to that in Figure~\figref[fig:model]). 
Before formally defining Markov counting systems, we introduce their key aspects. 
First, consider a population whose members are at any point in time in one (and only one) of $C$ possible stages (or compartments) of their lives, with stages belonging to finite collection $\mathcal{C}$.
Next, let the number of population members that are at stage $c$ at time $t$ define integer-valued random variables $X_c(t)$, which make up the system $\{\bm{X}(t)\} \equiv \{X_c(t): c \in \mathcal{C} \}$. 
Then, let the number of population members that have transitioned from stage $i$ to stage $j$ by time $t$ define non-decreasing, integer-valued random variables $N_{ij}(t)$, which in turn, for all pairs $\left(i,j\right)$ belonging to a collection of allowed transitions $\mathcal{T}$, define the collection of counting processes $\{\bm{N}(t)\} \equiv \{N_{ij}(t): (i,j) \in \mathcal{T}\}$. 
Next, let the collection \kwy{$\{\bm{N}(t)\}$} drive the dynamics of the system $\{\bm{X}(t)\}$ via the ``conservation of mass'' identity
\begin{eqnarray}
\label{eqn:mass}
X_{\MultInd}(t) 
= X_{\MultInd}(0)+\sum_{(i,\MultInd)\in\mathcal T} N_{i\MultInd}(t)-\sum_{(\MultInd,j) \in \mathcal T} N_{\MultInd j}(t), 
\end{eqnarray}
so that changes in $\{\bm{X}(t)\}$ are the result of changes in $\{\bm{N}(t)\}$. 
Mass conservation identity~\eqref{eqn:mass} restricts the transitions that can occur in $\{\bm{X}(t)\}$ as follows. 
Let $\N_0$ ($\N$) be the natural numbers including (excluding) zero and consider initial counts $\bm{n} \in \N^{\mathcal{T}}_0$ and initial system conditions $\bm{x} \in \N^{\mathcal{C}}_0$. 
For any given increments of the collection of counts $\bm{\ell} \equiv \; \{ \ell_{ij} : (i,j) \in \mathcal T \} \in \N^{\mathcal T}$, the system $\{\bm{X}(t)\}$ must make transitions $\bm{u} \equiv \{u_c : c\in\mathcal C\} \in \Z^{\mathcal C}$ with $u_c = \sum_{(i,c)\in \mathcal T} \ell_{ic} - \sum_{(c,j)\in \mathcal T} \ell_{cj}$. 
Finally, let the following \kwo{transition rates} define the Markov chain $\{\mathbf X(t),\mathbf N(t)\}$ 
\begin{eqnarray}
\label{eqn:MCS}
q(\bm{x}, \bm{\ell}) 
&\equiv& \lim\limits_{\dt \downarrow 0} \frac{P\Bigl(\bm{N}(t+\dt){=}\bm{n}+\bm{\ell}, \; \bm{X}(t+\dt){=}\bm{x}+\bm{u} \; | \bm{N}(t){=}\bm{n}, \bm{X}(t){=}\bm{x} \Bigr)}{\dt}.
\end{eqnarray}
Since the left hand side of \eqref{eqn:MCS} only depends on $\bm{x}$ (and not $\bm{n}$), $\{\bm{X}(t)\}$ is itself a continuous-time Markov chain and we call it a Markov counting system\footnote{The transition rates in \eqref{eqn:MCS} are time homogeneous, since its left hand side does not depend on $t$. 
This homogeneity adds clarity to the concepts, results and proofs but can readily be relaxed.}, 
which we illustrate with the following \kwdb{example}. 
%INTEPARAGRAPH SPACE

%INTEPARAGRAPH SPACE
\kwdb{Figure~\figref[fig:model]} defines a Markov counting system by relying on the concepts of \kwy{marginal transition rates} and of \kwo{pairwise transition rates}, which are necessary for its \kwb{interpretation} and which are also key to study the effect of \kwp{correlated external noise}. 
Consider the rate at which $k$ population members simultaneously undergo a transition of the $ij$-type (regardless of whether other members undergo other transitions), which can be defined as $q_{ij}(\bm{x}, k) \equiv \sum_{\bm{\ell}: \ell_{ij}= k}{q(\bm{x},\bm{\ell})}$ for $k \in \N$ and which we call the $(i,j)$ marginal transition rate. 
\kwy{Marginal transition rates} of size one $q_{ij}(\bm{x}, 1)$ are the labels on the arrows in Figure~\figref[fig:model]. 
Marginal rates of sizes greater than one do not appear in Figure~\figref[fig:model] because they are assumed to be zero. 
Another assumption needed to interpret Figure~\figref[fig:model] is that there are no co-jumps of different types. 
To formalize this second assumption, consider the rate at which $\bm{k} = (k_{ij},k_{i^{\prime}j^{\prime}})$ population members simultaneously undergo transitions of the $(i,j)$ and $(i^{\prime},j^{\prime})$ types (regardless of whether other members undergo other transitions), which can be defined as $q_{ij, i^{\prime}j^{\prime}}(\bm{x}, \bm{k}) \equiv \sum_{\bm{\ell}: \ell_{ij}= k_{ij}, \ell_{i^{\prime}j^{\prime}}= k_{i^{\prime}j^{\prime}}}{q(\bm{x},\bm{\ell})}$ for $\bm{k} \in \Big\{ \N_{0}^{2} - (0,0) \Big\}$ and which we call the $(i,j)-(i^{\prime},j^{\prime})$ pairwise transition rate. 
Requiring all \kwo{pairwise transition rates} to satisfy $q_{ij, i^{\prime}j^{\prime}}(\bm{x},(1,0)) = q_{ij}(\bm{x},1)$ and $q_{ij, i^{\prime}j^{\prime}}(\bm{x},(0,1)) = q_{i^{\prime}j^{\prime}}(\bm{x},1)$ guarantees no co-jumps and allows \kwb{interpreting} figures such as Figure~\figref[fig:model] as Markov chains \citep[see][]{anderson1991, bremaud1999, jacquez1996}. 
%For example, such interpretation implies that the rate at which an individual is first infected by strain 1, i.e., transitions from $S$ to $I_{1}$, is proportional to the rate at which an individual is infected by strain 1 after gaining some immunity by an infection from strain 2, i.e., transitions from $S_{1}$ to $I_{1}^{*}$. 
%This proportion is $q_{S_{1}I_{1}^{*}}(\bm{x}, 1) = (1-\gamma) q_{SI_{1}}(\bm{x}, 1) = (1-\gamma)  \lambda_{1}\xi_{1}$. 
%The interpretation also implies that such transitions happen one at a time, i.e., $q_{S_{1}I_{1}^{*}}(\bm{x}, k) =  q_{SI_{1}}(\bm{x}, k) = 0$ for $k > 1$, and never in groups of individuals, i.e., all $q_{SI_{1}, S_{1}I_{1}^{*}}(\bm{x}, (k_{SI_{1}}, k_{S_{1}I_{1}^{*}})) = 0$. 
Such interpretation also assumes that all non-zero rates $q_{ij}(\bm{x}, 1)$ are deterministic functions of the chain state $\bm{x}$ and not subject to \kwp{external noise}. 
%INTEPARAGRAPH SPACE
%
%INTEPARAGRAPH SPACE
%
%%%%%%%%%%%%%%%%%%%%%%%%%%%%%%%%%%%%%%%%%%%%%%%%%%%%
%%%%%%%%%%%%%%%%%%%%%%%%%%%%%%%%%%%%%%%%%%%%%%%%%%%%
%
% result section mark
\section{Correlated external noise in bivariate death Markov counting systems} 
\label{sec:bivar} 
Introducing \kwp{correlated noise} to the rates is easier if one considers two independent death processes and results in both co-jumps and in infinitesimal covariance as stated in \kwy{Proposition~\ref{pro:bivar}} below, which considers a common \kwo{multiplicative gamma} external noise, \kwb{a common death rate}, and which will later be useful when considering \kwdb{general Markov counting systems}. 
\kwy{Proposition~\ref{pro:bivar}} was proved in \cite{breto2011} and we state it here to make the paper self-contained. 
We state it in terms of our notation for Markov counting systems, after introducing our notation for the noise. 
The \kwo{noise} affecting the individual death rates is assumed to be continuous-time white noise obtained from a gamma process, which is also the choice of noise in \cite{breto2009} and \cite{shrestha2011}. 
Gamma white noise is defined as $\{\xi(t)\} \equiv \{d\Gamma(t)/dt\}$ with $\Gamma(t) \sim \text{Gamma}\,(t/\tau,\tau)$, $E[\Gamma(t)]=t$, and $V[\Gamma(t)]=\tau t$, so that $\tau$ parameterizes the magnitude of the noise.

\begin{proposition}[Proposition~7 of \citealp{breto2011}]\label{pro:bivar}
Consider the bivariate Markov counting system $\{ \bm{Y}(t) \} \equiv \big\{ \big( Y_1(t), Y_2(t) \big) \big\}$ defined by counting processes $\big \{ \big( N_{Y_{1}D}(t), N_{Y_{2}D}(t) \big) \big\}$ through mass conservation equations
\[ Y_{i}(t) = Y_{i}(0) - N_{Y_{i}D}(t),\]
and by transition rates $\q[y_{i},1][Y_{i}D][] = \delta y_{i} \Ind{0 < y_{i}}$
i.e., two independent linear death processes having equal individual death rate $\delta \in R^+$ and initial population sizes $Y_{i}(0)$. 
Consider subjecting both $\q[y_{i},1][Y_{i}D][]$ to a common gamma white noise, which defines counting processes $\{ N_{\tilde Y_{i} \tilde D}(t)\}=\{N_{Y_{i}D}\big(\Gamma(t)\big)\}$ and the corresponding Markov counting system $\{ \tilde{\bm{Y}}(t) \}$, i.e., two death processes each having stochastic rate $\delta\xi(t)$ . 
Then, the transition rates of $\{ \tilde{\bm{Y}}(t) \}$ correspond to 
\begin{comment}
marginal transition rates of univariate binomial gamma processes (as in \citealt{breto2011}), letting $k \leq x_{i}$, 
\begin{eqnarray}
\nonumber q_{X_{i}D} \big(\bm{x},k\big) &=& {x_{i} \choose k} \sum\limits_{j=0}^{k}{{k \choose j}}
(-1)^{k - j + 1}\tau^{-1} \ln{\bigl(1 + \delta \tau (x_i - j)\bigr)}
\end{eqnarray}
and 
\end{comment}
pairwise transition rates, letting $k_{i}=k_{\tilde Y_{i} \tilde D} \in \big\{ \N_{0}^{2} - (0,0): k_{i} \leq \tilde y_{i} \big\}$,
\begin{eqnarray}
\nonumber q_{ \tilde Y_{1} \tilde D, \tilde Y_{2} \tilde D} \big((\tilde y_1, \tilde y_2),(k_{1}, k_{2})\big) &=& {\tilde y_{1} \choose k_1} { \tilde y_{2} \choose k_2} \sum\limits_{j=0}^{k_1+k_2}{{k_1+k_2 \choose j}}
(-1)^{k_1+k_2 - j + 1}\tau^{-1} \ln{\bigl(1 + \delta \tau (\tilde y_1 + \tilde y_2 - j)\bigr)}.
\end{eqnarray}
Furthermore, the infinitesimal covariance between $\{ N_{\tilde Y_{1} \tilde D}(t)\}$ and $\{ N_{\tilde Y_{2} \tilde D}(t)\}$ is
\begin{eqnarray}
\nonumber
%\sigma_{d\tilde N_{X_{2}D}d\tilde N_{X_{2}D}}(x_{1},x_{2}) &=& 
\lim_{\dt \downarrow 0}\dt^{-1} Cov \Biggl[N_{\tilde Y_{1} \tilde D}(t+\dt) - N_{\tilde Y_{1} \tilde D}(t), \; N_{\tilde Y_{2} \tilde D}(t+\dt) - N_{\tilde Y_{2} \tilde D}(t) \; \bigg| \; \tilde{\bm{Y}}(t)=\tilde{\bm{y}} \Biggr] %\\
&=&
%\nonumber
\tilde y_{1} \tilde y_{2} \tau^{-1}\ln\Biggl(\frac{(1+\delta\tau)^2}{1+2\delta\tau}\Biggr)
%\hspace{3mm}
> %\hspace{2mm}
0.
\end{eqnarray}
\end{proposition}
%The proof of Proposition~\ref{pro:bivar} in \cite{breto2011} is based on the well-known result that $N_{X_{i}D}(t)$ follows a binomial distribution \citep{bharucha1960}. 
The \kwb{common rate assumption} of Proposition~\ref{pro:bivar} can be relaxed at the cost of more complex closed-form expressions for the covariance and for the pairwise rates. 
Although different death rates are assumed in the interacting death processes of Figure~\figref[fig:model], they will be assumed to be equal for the sake of simplicity when we illustrate in Section~\ref{sec:application} our results for \kwdb{general Markov counting systems} (which hold regardless of whether individual death rates are equal). 
%
%%%%%%%%%%%%%%%%%%%%%%%%%%%%%%%%%%%%%%%%%%%%%%%%%%%%
%%%%%%%%%%%%%%%%%%%%%%%%%%%%%%%%%%%%%%%%%%%%%%%%%%%%
%
% result section mark
\section{Correlated external noise in general Markov counting systems} 
\label{sec:inf-cov} 
Consider generalizing Proposition~\ref{pro:bivar} to \kwdb{general Markov counting systems}, which will lead us to defining \kwp{infinitesimal covariances} of such general systems. 
Consider a \kwdb{general Markov counting system} as defined by transition rates~\eqref{eqn:MCS} that satisfies the standard assumptions to interpret Figure~\figref[fig:model] of neither multiple jumps nor co-jumps and denote such system by $\{ \bm{W}(t) \}$. 
Consider now subjecting some (or all) transition rates of $\{ \bm{W}(t) \}$ to a collection of (possibly correlated, not necessarily gamma) white noises derived from $\{ \bm{ \Gamma }(t)\}$ (analogously to Section~\ref{sec:bivar}) and call the resulting process $\{ \tilde{\bm{W}} (t) \}$. 
The transitions rates of $\{ \tilde{\bm{W}} (t) \}$ could be obtained by integrating out the noises $\{\bm{\Gamma}(t)\}$ from the (now randomized) transition probabilities appearing inside the limit in~\eqref{eqn:MCS} as follows. 
Consider the collection of randomized time increments $\bm{H} \equiv \{H_{ij}: (i,j) \in \mathcal{T} \}$, necessarily with $E[H_{ij}]=\dt$. 
The nature of $H_{ij}$ depends on whether $q_{ij}(\bm{w},1)$ is subject to noise: if yes, then $H_{ij}$ is the corresponding noise random variable with density $f_{H_{ij}}$; if not, then it is the degenerate random variable $H_{ij}=\dt$. 
Then, provided they exist, the transition rates of $\{\tilde{\bm{W}}(t)\}$ are  
\begin{eqnarray}
\label{eqn:rates-correlated}
q(\tilde{\bm{w}}, \bm{\ell}) 
&=& \lim \limits_{\dt \downarrow 0} \int \frac{P\Bigl(\bm{N}(t+ \bm{s} ){=}\bm{n}+\bm{\ell}, \; \bm{W}(t+\bm{s}){=}\bm{w}+\bm{u} \; | \bm{N}(t){=}\bm{n}, \bm{W}(t){=}\bm{w} \Bigr)}{\bm{s}} f_{\bm{H}}(\bm{s})d\bm{s}.
\end{eqnarray}
While such direct integration of the noise in equation~\eqref{eqn:rates-correlated} was straightforward for the bivariate death process of Proposition~\ref{pro:bivar}, it is not so straightforward for more sophisticated models, like the one in Figure~\figref[fig:model] \citep{breto2009} or similar models \citep{shrestha2011}. 
Hence, instead of obtaining the new transition rates by direct integration, we propose constructing such new rates by directly specifying transition rates that produce the same infinitesimal covariance as that produced by introducing noise to appropriate bivariate systems.  
In the case of Figure~\figref[fig:model], noises $\xi_{i}(t)$ affect the rate of death processes $\{ N_{SI_{i}} \}$ and $\{ N_{S_{i} I_{i}^{*}} \}$. 
In this case, our proposal amounts to specifying a new set of transition rates for the system represented by Figure~\figref[fig:model] such that the infinitesimal covariance between $\{ N_{SI_{i}} \}$ and $\{ N_{S_{i} I_{i}^{*}} \}$ matches the covariance given by Proposition~\ref{pro:bivar} between $\{ N_{Y_{1}D} \}$ and $\{ N_{Y_{2}D} \}$. 
To do this, we first derive closed-form expressions for the \kwp{infinitesimal covariances} between two counting processes involved in a Markov counting system. 
%
%
%%%%%%%%%%%%%%%%%%%%%%%%%%%%%%%%%%%%%%%%%%%%%%%%%%%%
%%%%%%%%%%%%%%%%%%%%%%%%%%%%%%%%%%%%%%%%%%%%%%%%%%%%
%
% solution subsection mark
\subsection{Infinitesimal covariance of Markov counting systems} 
\label{sSec:inf-cov} 
Define the \kwp{infinitesimal covariances} of a general Markov counting system $\{\bm{X}(t)\}$ as defined in Section~\ref{sec:notation} as the collection
\[
\{\bm{\sigma}_{d\bm{X}}(\bm{x})\} \equiv \left\{ \sigma_{d\bm{X}}^{ij,i^{\prime}j^{\prime}}(\bm{x}): (i,j)\neq(i^{\prime},j^{\prime}) \in \mathcal{T}\right\}\] 
of infinitesimal covariances between counting processes $\{N_{ij}(t)\}$ and $\{N_{i^{\prime}j^{\prime}}(t)\}$:
\begin{eqnarray}
\label{eq:inf-cov}
\sigma_{d\bm{X}}^{ij,i^{\prime}j^{\prime}}(\bm{x})
& \equiv& 
\lim_{\dt \downarrow 0}\dt^{-1} Cov \Bigl[N_{ij}(t+\dt) - N_{ij}(t), \; N_{i^{\prime}j^{\prime}}(t+\dt) - N_{i^{\prime}j^{\prime}}(t) \; \big| \; \bm{X}(t)=\bm{x} \Bigr].
\end{eqnarray}
Our closed-form expressions below require one moment existence condition. 
Similar conditions were required by Theorem~1 of \cite{breto2011} to provide closed-form expressions for the infinitesimal mean and variance of Markov counting processes. 
Since \cite{breto2011} call their condition for the mean $(P1^{\star})$ and that for the variance $(P2^{\star})$, we shall call our condition for covariances $(P3^{\star})$.  
$(P3^{\star})$ is related to the number of transitions occurring in the Markov counting system over a time interval. 
This number of transitions is related not only to the sizes of the increments of each counting process $\{N_{ij}(t)\}$ but also to the overall rate at which these increments occur, which we call the rate function of the Markov counting system and define as 
\begin{eqnarray}
\label{eq:rate-fn}
\lambda(\bm{x}) 
\equiv \lim\limits_{\dt \downarrow 0} \frac{1 - P\left(N_{ij}\left(t+\dt\right) - N_{ij}\left(t\right) =0 \; \; \mbox{for all $\left(i,j\right) \in \mathcal{T} $}\; \big | \; \bm{X}\left(t\right)=\bm{x}\right)}{\dt}.
\end{eqnarray}
This quantity is also know as the intensity of the process in the point process literature \citep{daley2003}. 
If the rate function satisfies that $\lambda(\bm{x}) = \sum_{\bm{\ell}}{q(\bm{x},\bm{\ell})} < \infty$ for all $\bm{x}$, then the process is said to be stable and conservative. 
%INTEPARAGRAPH SPACE

%INTEPARAGRAPH SPACE
Consider stochastically bounding the rate function and the increment of each pair of counting processes over $[t,t+\bar h]$ by:
\begin{align}
\label{eq:bounds}
\bar\Lambda(t) 
&\equiv \sup_{t\le s\le t+\bar h}\lambda \Big( \bm{X}(s) \Big),
& \bar{Z}_{ij, i^{\prime} j^{\prime}}(t) &\equiv \sup{ \Big\{ \sup\limits_{t \leq s \leq t+\dtt}{dN_{ij}(s)}, \sup\limits_{t \leq s \leq t+\dtt}{dN_{i^{\prime}j^{\prime}}(s)} \Big \} }.
\end{align}
A combination of these two bounds gives the following property: 
\begin{itemize}
\item[\bf P3$^{\star}$.] For each $t$, $\bm{x}$ and $(i,j) \neq (i^{\prime}, j^{\prime})$ there is some $\bar h>0$ such that 
$E \Big[ \bar Z_{ij, i^{\prime}, j^{\prime}}^{2}(t)\bar \Lambda(t)|\bm{X}(t)=\bm{x} \Big] <\infty$. 
\end{itemize}
Property $(P3^{\star})$ requires that the Markov counting system does not have an explosive behaviour and holds, for example, for SIR-type compartmental models like that of Figure~\figref[fig:model], as shown in Section~\ref{sec:application}. 
It suffices to guarantee that infinitesimal covariances exist and that are given by the expression in Theorem~\ref{thm:inf-cov} below. 
%Supposing $(P3^{\star})$ formally extends the univariate theorem linking dispersion to the sample path by characterizing the infinitesimal dispersion of MCSs. 
%In the proof of the corollary we use Theorem~\ref{thm:mcs-inf-moments} which identifies the infinitesimal moments of a MCS. 
%
\begin{theorem}[Infinitesimal covariances of a Markov counting system]
\label{thm:inf-cov}
Let $\{\bm{X}(t)\}$ be a time homogeneous Markov counting system defined by counting processes $\{\bm{N}(t)\}$ and by transition rates $q(\bm{x}, \bm{\ell})$ as in~\eqref{eqn:MCS} that is stable and conservative. 
Supposing (P3$^\star$), the infinitesimal covariance between $\{N_{ij}(t)\}$ and $\{N_{i^{\prime}j^{\prime}}(t)\}$ is 
\[\sigma_{d\bm{X}}^{ij,i^{\prime}j^{\prime}}(\bm{x}) = \sum_{\bm{k}} k_{ij} k_{i^{\prime}j^{\prime}}  q_{ij, i^{\prime}j^{\prime}}(\bm{x},\bm{k}).\]
\end{theorem} 
\kwdb{Theorem~\ref{thm:inf-cov}} generalizes Theorem~1 of \cite{breto2011} to covariances and is proved in Appendix A. 
% ref2 comment 12
%Our proof uses the dominated convergence theorem and (P1$^{\star}$) and (P2$^{\star}$) require integrability of the dominating functions. 
%Then, the $\dt$ limit can be passed inside the moments in~\eqref{eq:dispersion}.
%
%
%INTEPARAGRAPH SPACE
%
%INTEPARAGRAPH SPACE
%
%
%%%%%%%%%%%%%%%%%%%%%%%%%%%%%%%%%%%%%%%%%%%%%%%%%%%%
%%%%%%%%%%%%%%%%%%%%%%%%%%%%%%%%%%%%%%%%%%%%%%%%%%%%
%
% application section mark
\section{Transition rates of SIR-type models subject to external correlated noises} 
\label{sec:application} 
\kwdb{Theorem~\ref{thm:inf-cov}} can be used to show that, after minimal \kwy{simplifications} to add clarity to our contribution, the system represented in Figure~\figref[fig:model] can be defined by \kwo{transition rates} that \kwb{reproduce the effects} (identified in Proposition~\ref{pro:bivar}) of correlated noises. %, as well as to \kwb{interpret such noises} in terms of infinitesimal covariances in deterministic environments (instead of in terms transition rates in random environments). 
First, consider the Markov chain $\{ \bm{Z}(t) \} \equiv \big\{ \big( S(t), I_{1}(t), I_{2}(t), S_{1}(t), S_{2}(t), I_{1}^{*}(t), I_{2}^{*}(t), R(t) \big) \big\}$ specified by Figure~\figref[fig:model]. 
The standard interpretation of Figure~\figref[fig:model] gives the transition rates for $\{ \bm{Z}(t) \}$ in Table~1 (e.g., letting the $\xi_{i}$ be deterministic constants).  
Next, before considering adding noise to $\{ \bm{Z}(t) \}$, we make the following simplifications to $\{ \bm{Z}(t) \}$ so that the contributions in this paper can be presented more clearly: 
(i) instead of a time-inhomogeneous birth rate $b(t)$, births compensate deaths so that the total population size remains constant and is equal to $P < \infty$, as in \cite{shrestha2011}; 
(ii) instead of a time-inhomogeneous infection rate $\beta(t)$ within $\lambda_{i}$, this rate is constant and equal to $\beta$; and 
(iii) instead of different transition rates from $S$ to $I_{i}$ and from $S_{i}$ to $I_{i}^{*}$, these rates are both $\lambda_{i} \xi_{i}(t)$, i.e., $\gamma=0$.
\kwy{Simplifications} (i) and (ii) impose time-homogeneity, which allows for a simpler notation in the rest of the paper. 
Simplification (iii) allows for simpler transition rates and covariance closed-form expressions. 
Now, let stable, conservative Markov chain $\{ \tilde{\bm{Z}}(t)\}$ be defined by the rates of Table~1 modified according to (i)--(iii) and by non-zero pairwise rates of transitions involving the $\xi_{i}$ equal to those in Proposition~\ref{pro:bivar} as follows: 
\begin{eqnarray}
\nonumber
q_{\tilde S \tilde I_{i}, \tilde S_{i} \tilde I_{i}^{*}}( \tilde{\bm{z}},\bm{k}) 
&=& 
{ \tilde s \choose k_1} {\tilde s_{i} \choose k_2} \sum\limits_{j=0}^{k_1+k_2}{{k_1+k_2 \choose j}}
(-1)^{k_1+k_2 - j + 1}\tau^{-1} \ln{\bigl(1 + \lambda_{i} \tau (\tilde s + \tilde s_i - j)\bigr)}. 
\end{eqnarray}
System $\{ \tilde{\bm{Z}}(t)\}$ defined by \kwo{these rates} satisfies $(P3^{\star})$, since letting the fixed population size be $\tilde P$,
\begin{eqnarray}
\nonumber
\lambda \Big( \tilde{\bm{Z}}(t) \Big)
 &=& m \big( \tilde P - \tilde S(t) \big) + r \Big( \tilde I_{1}(t) + \tilde I_{2}(t) + \tilde I_{1}^{*}(t) + \tilde I_{2}^{*}(t)  \Big) + 
% \sum\limits_{i=1,2} \Biggl[ \sum\limits_k \q[\bm{Z}(t),k][SI_{i}][] + \sum\limits_k \q[\bm{Z}(t),k][S_{i}I_{i}^{*}][] + \sum\limits_{\bm{k}} \q[\bm{Z}(t),k][SI_{i}, S_{i}I_{i}^{*}][]  \Biggr] \\
\sum\limits_{\bm{k}} \q[\tilde{\bm{Z}}(t),\bm{k}][\tilde S \tilde I_{1}, \tilde S_{1} \tilde I_{1}^{*}][] + \sum\limits_{\bm{k}} \q[\tilde{\bm{Z}}(t),\bm{k}][\tilde S \tilde I_{2}, \tilde S_{2} \tilde I_{2}^{*}][]\\
\label{eqn:rate-bound}
&\leq& \Bigg(m + r  + \lambda_{1} +  \lambda_{2} \Bigg) \tilde P
%
%\lambda_{1} \big( S(t) + S_{1} \big)  + \lambda_{2}  \big( S(t) + S_{2} \big),
%
\end{eqnarray}
where the inequality follows by substituting all compartments by $\tilde P$ and because 
\[ 
\sum_{\bm{k}} \q[\tilde{\bm{Z}}(t),\bm{k}][\tilde S \tilde I_{i}, \tilde S_{i} \tilde I_{i}^{*}][] = \tau^{-1} \ln{\Big(1 + \tau \lambda_{i} \big( \tilde S(t) + \tilde S_{i} \big) \Big)} \leq \lambda_{i} \tilde P 
\] 
(as follows from the properties of the binomial gamma process of \citealp{breto2011}). 
Since \eqref{eqn:rate-bound} is not time-varying, it also bounds $\bar \Lambda (t)$ involved in $(P3^{\star})$. 
Similarly, $\tilde P$ is an upper bound for the increments: $\bar Z_{ij, i^{\prime} j^{\prime}}(t) \leq \tilde P$, so that 
\bean
E \left[ \bar Z_{ij, i^{\prime} j^{\prime}}^2 (t) \bar \Lambda(t) \; | \; \tilde{\bm{Z}}(t) = \tilde{\bm{z}} \right] 
\;\; \leq \;\;
\Bigg(m + r  + \lambda_{1} +  \lambda_{2} \Bigg) \tilde P^{3} \eean

$(P3^{\star})$ can be analogously verified for the bivariate system $\{ \tilde{\bm{Y}}(t)\}$ of Proposition~\ref{pro:bivar} by assuming, for example, that the initial population sizes are deterministic, i.e., for fixed $y_{i}(0) = \tilde y_{i}(0)$ 
\[ 
\lambda \Big( \tilde{\bm{Y}}(t) \Big) = \sum_{\bm{k}} \q[\tilde{\bm{Y}}(t),\bm{k}][\tilde Y_{1} \tilde D, \tilde Y_{2} \tilde D][] = \tau^{-1} \ln{\Big(1 + \tau \delta \big( \tilde Y_{1}(t) + \tilde Y_{2}(t) \big) \Big)} \leq \delta \big( \tilde  y_{1}(0) + \tilde y_{2}(0) \big).
\] 
Hence, it follows directly from Theorem~\ref{thm:inf-cov} that 
\begin{eqnarray}
\label{eqn:fx}
\sigma_{d\bm{ \tilde Z}}^{\tilde S \tilde I_{i}, \tilde S_{i} \tilde I_{i}^{*}}(\tilde{\bm{z}}) 
&=&
\tilde s \tilde s_{i} \tau^{-1}\ln \Biggl( \frac{(1+\lambda_{i} \tau)^2}{1+2 \lambda_{i} \tau} \Biggr) > 0.
\end{eqnarray}
Equation~\eqref{eqn:fx} can be interpreted as follows. 
First, the effects of introducing correlated noises in the bivariate system of Proposition~\ref{pro:bivar} can be reproduced in more general systems. 
Second, it permits an alternative interpretation of correlated noises $\xi_{i}(t)$ in a non-random environment context. 
These noises have effectively been integrated out in $\{ \tilde{\bm{Z}}(t)\}$, which can be seen as a regular Markov chain in a non-random environment, with the caveat that it now allows for simultaneous co-transitions that drive the new infinitesimal correlations. 
\renewcommand{\arraystretch}{1.5}
\begin{table}[h]
\label{tab:simple}
\caption{Transition rates according the standard interpretation of Figure~\figref[fig:model] as a continuous-time Markov chain with rate function $\lambda_{\bm{Z}}(\bm{z}) \equiv \sum_{(i,j) \in \mathcal{T}}  q_{ij}(\bm{z},1) $ and with all marginal rates $q_{ij}(\bm{z},k)$ for $k>1$ and all pairwise transition rates $q_{ij, i^{\prime}j^{\prime}}(\bm{z},\bm{k})$ assumed to be zero.}
\centering
\begin{tabular}{c | l*{11}{c}}
\hline
\hline
$N_{ij}$ & $N_{SI_{i}}$ & $N_{I_{i}S_{i}}$ & $N_{S_{i}I_{i}^{*}}$ & $N_{I_{i}^{*}R}$ & $N_{BS}$ & $N_{SD}$ & $N_{I_{i}D}$ & $N_{S_{i}D}$ & $N_{I_{i}^{*}D}$ & $N_{RD}$\\
\hline
$q_{ij}(\bm{z},1)$ & $\lambda_{i} \xi_{i}$ & $r$ & $(1-\gamma) \lambda_{i} \xi_{i}$ & $r$ & $b(t)$ & $m$ & $m$ & $m$ & $m$ & $m$\\
\hline
\end{tabular}
\end{table}
%
%INTEPARAGRAPH SPACE
%
%INTEPARAGRAPH SPACE
%
%%%%%%%%%%%%%%%%%%%%%%%%%%%%%%%%%%%%%%%%%%%%%%%%%%%%
%%%%%%%%%%%%%%%%%%%%%%%%%%%%%%%%%%%%%%%%%%%%%%%%%%%%
%
¥% Acknowledgement section mark
\section*{Acknowledgements}
This work was supported by Spanish Government Project ECO2012-32401 and Spanish Program \emph{Juan de la Cierva} (JCI-2010-06898). 
%\end{comment}
%I thank the editor, associate editor and two anonymous referees for their suggestions which helped improve the paper. %led to substantial improvements in the paper. 
%I thank %Edward Ionides for insightful comments and 
%the editor and an anonymous referee for suggestions that helped improve the paper.
%
%% References
%%   \citep{key}         ==>>  [#]
%\bibliographystyle{apalike}
\bibliographystyle{elsarticle-harv}
%\bibliographystyle{elsarticle-num}
%\bibliographystyle{pjgsm}
%\bibliography{sctmcp-pois-bib,}
\bibliography{references}
%\bibliography{references}
%% \bibitem{key}...
% \end{thebibliography}
%
%
%
%% Appendix
\appendix
\renewcommand{\thesection}{Appendix \Alph{section}}
\renewcommand{\theequation}{\Alph{section}.\arabic{equation}}
\setcounter{equation}{0}
\setcounter{section}{0}
%
%%%%%%%%%%%%%%%%%%%%%%%%%%%%%%%%%%%%%%%%%%%%%%%%%%%%
%%%%%%%%%%%%%%%%%%%%%%%%%%%%%%%%%%%%%%%%%%%%%%%%%%%%
%
¥% appendix mark
\section{Proof of Theorem~\ref{thm:inf-cov}}
\begin{proof}
First, we prove that the infinitesimal covariance equals the infinitesimal cross-product under condition $(P3^{\star})$. 
Let $\Delta N_{ij}(t) \equiv N_{ij}(t + \dt) - N_{ij}(t)$ and analogously for all other subindices and counting processes in this proof. 
Since $(P3^{\star})$ implies $(P1^{\star})$ in Theorem~1 of \cite{breto2012}, it follows that $E \Bigl[ \Delta N_{ij}(t)\; \big| \; \bm{X}(t)=\bm{x} \Bigr] = \dt \sum\limits_{k} k q_{ij}(\bm{x},k) + o(\dt)$. 
Then
\begin{eqnarray}
\nonumber
\lim_{\dt \downarrow 0}\dt^{-1} Cov \Bigl[ \Delta N_{ij}(t), \; \Delta N_{i^{\prime}j^{\prime}}(t) \; \big| \; \bm{X}(t)=\bm{x} \Bigr] 
&=& \lim_{\dt \downarrow 0}\dt^{-1} \Bigg\{ E \Bigl[ \Delta N_{ij}(t) \Delta N_{i^{\prime}j^{\prime}}(t) \; \big| \; \bm{X}(t)=\bm{x} \Bigr] - o(\dt) \Bigg\}\\
%
%\nonumber
%&& - E \Bigl[ \Delta N_{ij}(t)\; \big| \; \bm{X}(t)=\bm{x} \Bigr] E \Bigl[ \Delta N_{i^{\prime}j^{\prime}}(t)\; \big| \; \bm{X}(t)=\bm{x} \Bigr] \Bigg\}\\ 
%
\nonumber
&=& \lim_{\dt \downarrow 0}\dt^{-1} E \Bigl[ \Delta N_{ij}(t) \Delta N_{i^{\prime}j^{\prime}}(t) \; \big| \; \bm{X}(t)=\bm{x} \Bigr].
\end{eqnarray}
The rest of this proof follows closely the proof of Theorem~1 in \cite{breto2012}. 
While that theorem provided expressions for the infinitesimal mean and variance of $\{\bm{X}(t)\}$, this one provides them for the covariances. 
%INTEPARAGRAPH SPACE

%INTEPARAGRAPH SPACE
All probabilities and expectations in this proof are conditional on $\bm{X}(t) = \bm{x}$ (in addition to other conditioning, where appropriate). 
Define the following: (i)~let $\{\bar{N}_{ij}(t)\}$ be a process such that, conditional on $\bar{\Lambda}(t)$ and $\bar{Z}(t) \equiv \bar{Z}_{ij, i^{\prime} j^{\prime}}(t)$, realizations of $\{\bar{N}_{ij}(t)\}$ are those of a compound Poisson process \citep{cox1980} with Poisson event rate $\bar{\Lambda}(t)$ and degenerate jump or batch size distribution \citep{daley2003} with mass one at $\bar{Z}(t)$, i.e., a process with jumps arriving according to the Poisson process and for which the size of the jumps is $\bar{Z}(t)$; and 
(ii)~let $S$ be the event that there is exactly one transition time occurring in the interval $[t,t+\dt]$ in the Markov counting system $\{\bm{X}(t)\}$. 
%At this single transition time, different simultaneous transitions could increase more than one of the counting processes $\{N_{ij}(t)\}$. 
%Also, these increments could be of size greater than one if the processes are compound. 
Then,
\begin{eqnarray}
\label{eq:s-and-sc}
E \left[ \Delta N_{ij}(t) \Delta N_{i^{\prime}j^{\prime}}(t) \right] 
&=& E \left[ \Delta N_{ij}(t) \Delta N_{i^{\prime}j^{\prime}}(t) \Ind{S} \right] + E \left[ \Delta N_{ij}(t) \Delta N_{i^{\prime}j^{\prime}}(t) \Ind{S^c} \right].
\end{eqnarray} 
Consider the first term on the right hand side of~\eqref{eq:s-and-sc}. 
Let $S_{ij, i^{\prime}j^{\prime}} \subset S$ be the event that there is exactly one transition time occurring in the interval $[t,t+\dt]$ in the MCS $\{\bm{X}(t)\}$ and that this transition increases both the $\{N_{ij}(t)\}$ and the $\{N_{ i^{\prime}j^{\prime}}(t)\}$ processes (and possibly other processes). 
Then, in~\eqref{eq:s-and-sc}, letting $\bm{k} \equiv (k_{ij}, k_{i^{\prime}j^{\prime}})$ 
\begin{eqnarray}
\nonumber
E[\Delta N_{ij}(t) \Delta N_{i^{\prime}j^{\prime}}(t) \Ind{S}] 
&=& E[\Delta N_{ij}(t) \Delta N_{i^{\prime}j^{\prime}}(t) | S_{ij, i^{\prime}j^{\prime}}] \;
\times \; \; \, P(S_{ij, i^{\prime}j^{\prime}}|S) \, \; 
\times \; \; \, P(S)\\
\label{eq:ross-result}
&=& \sum\limits_{\bm{k}} k_{ij} k_{i^{\prime}j^{\prime}} \frac{\q[\bm{x},\bm{k}][ij, i^{\prime}j^{\prime}][]}{\sum\limits_{\bm{k}} \q[\bm{x},\bm{k}][ij, i^{\prime}j^{\prime}][]} 
\;\; \times \frac{\sum\limits_{\bm{k}} \q[\bm{x}, \bm{k}][ij, i^{\prime}j^{\prime}][]}{\lambda(\bm{x})} 
\times \Bigl[\dt \lambda(\bm{x}) + o(\dt)\Bigr]\\ 
&=& \dt \sum\limits_{\bm{k}} k_{ij} k_{i^{\prime}j^{\prime}}  \q[\bm{x},\bm{k}][ij, i^{\prime}j^{\prime}][] + o(\dt)
\end{eqnarray}
where $P(S)$ in~\eqref{eq:ross-result} follows by a standard result on Markov chains \citep[see for example:][page 492]{ross1996}.
%INTEPARAGRAPH SPACE

%INTEPARAGRAPH SPACE
To finish the proof, we show that the second term on the right hand side of~\eqref{eq:s-and-sc} disappears infinitesimally. 
Let $\bar S$ be the event that there is exactly one transition time occurring in the interval $[t,t+\dt]$ in the compound Poisson process $\{\bar N_{ij}(t)\}$. 
Since the random variable $\Delta N_{ij}(t) \Delta N_{i^{\prime}j^{\prime}}(t)$ is stochastically smaller than $\Big( \Delta \bar{N}_{ij}(t) \Big)^2$, 
\begin{eqnarray}
\nonumber
E \left[ \Delta N_{ij}(t) \Delta N_{i^{\prime}j^{\prime}}(t) \Ind{S^c} \right] 
&\leq& E\left[ \Big( \Delta \bar{N}_{ij}(t) \Big)^2 \Ind{\bar S^c} \right]\\ 
\label{eq:expected-value-indicator}
&=& E \Bigg[ E \Big[  \Big( \Delta \bar{N}_{ij}(t) \Big)^2  \big | \; \bar{\Lambda}(t),\bar{Z}(t) \Big] \Bigg] - E \Bigg[ E \Big[ \Big( \Delta \bar{N}_{ij}(t) \Big)^2  \Ind{\bar S} \; \big | \; \bar{\Lambda}(t),\bar{Z}(t) \Big] \Bigg]\\ 
\label{eq:poisson-prop-mean}
&=& \Bigg( E \Big[ \bar{Z}^2(t) \bar{\Lambda}(t) \dt \Big] + \underbrace{E \Big[ \bar{Z}^2(t) \bar{\Lambda}^2(t) \Big] \dt^2}_{=o(\dt)} \Bigg) - E \Big[ \bar{Z}^2(t) \bar{\Lambda}(t) \dt \exp{\{-\dt \bar{\Lambda}(t)\}} \Big]\\
\nonumber
&=& E \bigg[ \bar{Z}^2(t) \bar{\Lambda}(t) \dt \Big(1 - \exp \big\{-\dt \bar{\Lambda}(t) \big\} \Big) \bigg] + o(\dt) 
\end{eqnarray}
where~\eqref{eq:expected-value-indicator} follows as in~\eqref{eq:s-and-sc}, and~\eqref{eq:poisson-prop-mean} follows by the properties of the compound Poisson distribution. 
Since $\bar{z}^{2}\bar{\lambda}\Big(1 - \exp\big\{-\dt \bar{\lambda}\big\}\Big) \leq \bar{z}^{2}\bar{\lambda}$ and $E \left[ \bar{Z}^{2}(t) \bar{\Lambda}(t) \right]$ is assumed finite (note that the distribution of $\bar Z^{2}(t) \bar \Lambda(t)$ depends on $\dtt$ and not $\dt$), it follows by dominated convergence that
\begin{eqnarray}
\nonumber \lim\limits_{\dt \downarrow 0} \frac{E\bigg[ \bar{Z}^{2}(t) \bar{\Lambda}(t) \dt \Bigl(1 - \exp\big\{-\dt \bar{\Lambda}(t)\big\}\Bigr)\bigg]}{\dt} 
&=& E\biggl[\lim\limits_{\dt \downarrow 0} \bar{Z}^{2}(t) \bar{\Lambda}(t)\Bigl(1 - \exp\big\{-\dt \bar{\Lambda}(t)\big\}\Bigr)\biggr] = 0.
\end{eqnarray} 
\end{proof}
\end{document}